\documentclass[11pt,a4paper]{amsart}
\usepackage[english]{babel}
\usepackage{hyperref}
\usepackage{german}
\usepackage{amsmath}
\usepackage{hyperref}
\usepackage{amssymb}
\usepackage{amsthm}
\usepackage{fontenc}
\usepackage{graphicx}
\usepackage{enumerate}
\usepackage{a4wide}
\usepackage[arrow, matrix, curve]{xy}
\usepackage{color}
\numberwithin{equation}{section}

\theoremstyle{plain}
\newtheorem{satz}{Theorem}[section]
\newtheorem{lem}[satz]{Lemma}
\newtheorem{fol}[satz]{Corollary}
\newtheorem{prop}[satz]{Proposition}
\newtheorem{con}[satz]{Conjecture}

\theoremstyle{definition}
\newtheorem{defi}[satz]{Definition} 
\newtheorem{bsp}[satz]{Example}

\newtheorem*{Frage}{Question}

\theoremstyle{remark}
\newtheorem{nota}[satz]{Notation}
\newtheorem{bem}[satz]{Remark}

\newcommand{\maxk}[1]{\left\{#1\right\}}
\newcommand{\erz}[1]{\langle#1\rangle}
\newcommand{\La}{\lambda}
\newcommand{\LA}{\Lambda}
\newcommand{\n}{\nu}

\DeclareMathOperator{\reg}{reg}
\DeclareMathOperator{\codim}{codim}

\DeclareMathOperator{\height}{ht}

\begin{document}

\title[A combinatorial proof of the E-G conjecture for monomial curves]{A combinatorial proof of the Eisenbud-Goto\\ conjecture for monomial curves and\\ some simplicial semigroup rings}
\author{Max Joachim Nitsche}
\address{Max\mbox{\;}Planck\mbox{\;}Institute\mbox{\;}for\mbox{\;}Mathematics\mbox{\;}in\mbox{\;}the\mbox{\;}Sciences,\mbox{\;}Inselstrasse\mbox{\;}22,\mbox{\;}04103\mbox{\;}Leipzig,\mbox{\;}Germany}
\email{nitsche@mis.mpg.de}
\thanks{}
\date{\today}


\begin{abstract}

We will give a pure combinatorial proof of the Eisenbud-Goto conjecture for arbitrary monomial curves. Moreover, we will show that the conjecture holds for certain simplicial affine semigroup rings.

\end{abstract}

\maketitle

\section{Introduction}

In \cite{KMNSNAX} we established new bounds for the Castelnuovo-Mumford regularity of (homogeneous) seminormal simplicial affine semigroup rings, in particular, we showed that the Eisenbud-Goto conjecture \cite{KEG} holds in this case, that is, its Castelnuovo-Mumford regularity is bounded by its degree minus its codimension. One of the main keys to prove this was an idea of Hoa and St\"uckrad, namely, the regularity of simplicial affine semigroup rings can be computed in terms of the regularity of certain monomial ideals. We will again use this idea to give a pure combinatorial proof of the Eisenbud-Goto conjecture for arbitrary monomial curves and certain simplicial affine semigroup rings.\\

Let $K$ be a field, and let $R=K[x_1,\ldots,x_n]$ be a standard graded polynomial ring, that is, all variables $x_i$ have degree $1$. We define the \emph{Castelnuovo-Mumford regularity} (or \emph{regularity} for short) $\reg M$ of a finitely generated graded $R$-module $M$ by
$$
\reg M:=\max\maxk{i+a(H^i_{R_+}(M))\mid i\geq0},
$$
where $H^i_{{R}_+}(M)$ denotes the $i$-th local cohomology module of $M$ with respect to the homogeneous maximal ideal $R_+$ of $R$, and $a(H^i_{R_+}(M)):=\max\maxk{r\mid H^i_{R_+}(M)_r\not=0}$ with the convention $a(0)=-\infty$. Bounding the Castelnuovo-Mumford regularity of $M$ is interesting for various reasons, for example, the $i$-th syzygy module of $M$ can be generated by elements of degree smaller or equal to $\reg M+i$ by \cite{KEG}, moreover, the Hilbert function of $M$ agrees with its Hilbert polynomial for inputs bigger than $\reg M$ (see, for example, \cite{KTGOS}). In addition to this, one can even bound under certain assumptions the degrees in a minimal Gr\"obner bases of a homogeneous ideal $I$ of $R$ by $\reg R/I+1$ by a result of Bayer and Stillman \cite{KBSCDM}.\\

In the following we will consider homogeneous simplicial affine semigroup rings. We define the set $M_{d,\alpha}:=\{(u_{1},\ldots,u_{d})\in\mathbb N^d \mid \sum_{i=1}^{d}u_{i}=\alpha\}$ where $d,\alpha\in\mathbb N^+$. Let $B$ be a homogeneous simplicial affine semigroup, that is, up to isomorphism (see \cite{KMNSNAX}) we may assume that $B$ is the submonoid of $(\mathbb N^d,+)$ which is generated by a set $\{e_1,\ldots,e_d,a_1,\ldots,a_c\}\subseteq M_{d,\alpha}$, where
$$
e_1:=(\alpha,0,\ldots,0), e_2:=(0,\alpha,0,\ldots,0),\ldots,e_d:=(0,\ldots,0,\alpha).
$$
We may also assume that the integers $a_{i[j]},\:i=1,\ldots,c,\:j=1,\ldots,d$ are relatively prime, where $a_i=({a_i}_{[1]},\ldots,{a_i}_{[d]})$, moreover, we assume that $c\geq1$. By $K[B]$ we denote the affine semigroup ring associated to $B$. As usual we identify the ring $K[B]$ with the subring of the polynomial ring $K[t_1,\ldots,t_d]$ generated by monomials $t^b:=t_1^{b_1}\cdot\ldots\cdot t_d^{b_d}$ for $b=(b_1,\ldots,b_d)\in B$. We will always consider the positive $\mathbb Z$-grading on $K[B]$ which is induced by $\deg t^b=(\sum_{i=1}^d b_i)/\alpha$. Denote by $\deg K[B]$ the degree of $K[B]$ and by $\codim K[B]:=\dim_K K[B]_1-\dim K[B]=c$ its codimension; note that $\dim K[B]=d$. By $\reg K[B]$ we mean the regularity of $K[B]$ with respect to the canonical $R$-module structure which is given by the homogeneous surjective $K$-algebra homomorphism
$$
\pi: R=K[x_1,\ldots,x_{d+c}]\twoheadrightarrow K[B],
$$
where $x_i\mapsto t_i^\alpha$ for $i=1,\ldots,d$ and $x_{d+j}\mapsto t^{a_j}$ for $j=1,\ldots,c$. Thus, $R/\ker\pi\cong K[B]$, where $\ker\pi$ is a homogeneous prime ideal. For (homogeneous) simplicial affine semigroup rings the Eisenbud-Goto conjecture comes down to the

\enlargethispage{0.14cm}

\begin{Frage}[Eisenbud-Goto {\cite{KEG}}]

Does $\reg K[B]\leq \deg K[B] - \codim K[B]$ hold?

\end{Frage}

Recall that from the general known results (see \cite{KMNSNAX}) one can deduce that this question has a positive answer if $\dim K[B]=2$ by Gruson, Lazarsfeld, and Peskine \cite{KGLP}; if $K[B]$ is Cohen-Macaulay or Buchsbaum by Treger, St\"uckrad and Vogel \cite{KEGCM,KEGBB}; and if $\deg K[B]\leq c+2$ by Hoa, St\"uckrad, and Vogel \cite{KHSVC2}. Moreover, the question has a positive answer if $c=2$ by Peeva and Sturmfels \cite{Kcodim2}; if $c\leq{\rm deg}K[B]/\alpha$ by Hoa and St\"uckrad \cite{KHSCM}; if ${\rm deg}K[B]=\alpha^{d-1}$ and $\alpha\leq d-1$ again by \cite{KHSCM}; if $K[B]$ has an isolated singularity, equivalently, $B$ contains all elements of $M_{d,\alpha}$ such that one coordinate is equal to $\alpha-1$, that is, it contains all elements of type $(0,\ldots,\alpha-1,\ldots,1,\ldots,0)$, by Herzog and Hibi \cite{KCMHH}; and if $K[B]$ is seminormal by \cite{KMNSNAX}. We also refer to the paper of Lazarsfeld \cite{KLSMSD3} for a proof of the Eisenbud-Goto conjecture for smooth surfaces in characteristic zero, to Ran \cite{KRLDG} for a proof of the conjecture for certain smooth threefolds in characteristic zero, and to Giaimo \cite{KDGEGCC} for a proof for connected reduced curves. Thus, the Eisenbud-Goto conjecture is widely open, even if we restrict ourself to simplicial affine semigroup rings.\\

Every homogeneous affine semigroup such that its semigroup ring has dimension $2$ is simplicial and is therefore isomorphic to some $B$ with $\dim K[B]=2$. The ring $K[B]$ is isomorphic to the coordinate ring of a (projective) monomial curve of degree $\alpha$ in $\mathbb P^{c+1}$ if $\dim K[B]=2$, thus, for simplicity we will also use the term monomial curve instead of $K[B]$ in this situation. So far every proof of the Eisenbud-Goto conjecture for monomial curves uses some techniques of the Gruson, Lazarsfeld, and Peskine proof. In contrast to this approach Bruns, Gubeladze, and Trung asked in \cite{KBGTPAS} for a combinatorial proof of the conjecture in this case. We will completely answer this question by giving the first combinatorial proof of the Eisenbud-Goto conjecture for arbitrary monomial curves in Theorem~\ref{Kegmonomial}. Our proof is elementary and uses an idea of Hoa and St\"uckrad, namely, one can decompose $K[B]$ into a direct sum of certain monomial ideals, and hence compute its regularity in terms of the regularity of the monomial ideals. This becomes even more powerful in case of monomial curves, since all monomial ideals are contained in a polynomial ring in two variables - this enables us to read off the regularity of the monomial ideals by ordering its minimal monomial generators. Moreover, we will use the idea of sequences with $*$-property which has been developed in \cite{KMNSNAX} to give a proof of the conjecture in the seminormal case, roughly speaking, these sequences are useful to control degree of $K[B]$ in a certain way. This will further lead to a proof of the Eisenbud-Goto conjecture in case that one of the monomial ideals in the decomposition which determines the regularity of $K[B]$ is generated by at most two elements in Theorem~\ref{Keg2}; note that this result holds in any dimension.\\

By using the Gruson, Lazarsfeld, Peskine techniques L'vovsky \cite{KLV} showed that the regularity of a monomial curve is bounded by $ \#L+\#L'+1$, where $L$ and $L'$ are the longest and the second longest gap of $B$, this bound will be called the L'vovsky bound; a gap is a maximal set of consecutive integer points on the line $[(\alpha,0),(0,\alpha)]$ not belonging to $B$. Note that the L'vovsky bound is better than the Eisenbud-Goto bound, since $\deg K[B]-c=\sum (\#L) +1$ where the sum is taken over all gaps $L$ of $B$. Under the assumption $(1,\alpha-1),(\alpha-1,1)\in B$ we even get a better bound, namely, $\reg K[B]$ is bounded by $\#L+1$ where $L$ is the longest gap of $B$ by a result of Hellus, Hoa, and St"uckrad \cite{KCMRMC}. But even this combinatorial bound in \cite{KCMRMC} is far from sharp for $c\geq4$, see, for example, \cite[Introduction]{KDAMNP}. For other combinatorial results for monomial curves we refer to the paper of Thomas \cite{KHTBDMC} in which it is shown that the corresponding toric ideal can be generated by elements of degree smaller or equal to the L'vovsky bound plus $1$, and to the paper of Crupi and Utano \cite{KECEGMC} in which the Eisenbud-Goto bound has been proved for very special classes of monomial curves.\\

In Section~\ref{Kbasics} we will again recall the construction of the decomposition of the ring $K[B]$, moreover, we will develop the main tools which are needed to prove the assertions in Section~\ref{Kcase2} and in Section~\ref{Kmonomial}. For unspecified notation we refer to \cite{KBGPRKT, KEB}.

\enlargethispage{0.12cm}

\section{Basics}\label{Kbasics}

We define $A:=\erz{e_1,\ldots,e_d}$ to be the submonoid of $B$ generated by $e_1,\ldots,e_d$.
Denote by $G(A)$ and $G(B)$ the groups generated by $A$ and $B$ respectively. For $x\in \mathbb Z^d$ we denote by $x_{[i]}$ the $i$-th component of $x$, and for $x\in G(B)$ we define $\deg x:=(\sum_{i=1}^dx_{[i]})/\alpha$. We set 
$$
B_A:=\{x\in B\mid x-a\notin B~\forall a\in A\setminus\{0\}\};
$$ 
note that the set $B_A$ is always finite. Moreover, if $x\notin B_A$, then $x+y\notin B_A$ for all $x,y\in B$. We define the equivalence relation $\sim$ on $G(B)$ by $x\sim y$ if $x - y\in G(A)= \alpha\mathbb Z^d$. One can show that there are exactly $f:=\#(G(B)\cap D)$ equivalence classes on $G(B)$, $B$, and on $B_A$, where $D:=\{(x_{[1]},\ldots,x_{[d]})\in\mathbb Q^d\mid 0\leq x_{[i]}<\alpha~\forall i\}$. By $\Gamma_1,\ldots,\Gamma_f$ we denote the equivalence classes on $B_{A}$. For $t=1,\ldots,f$ we define
$$
h_t:=(\min\maxk{m_{[1]}\mid m\in \Gamma_t}, \min\maxk{m_{[2]}\mid m\in \Gamma_t},\ldots,\min\maxk{m_{[d]}\mid m\in \Gamma_t}).
$$
Note that by construction $x-h_t\in A$ for all $x\in\Gamma_t$, and hence $h_t\in G(B)\cap \mathbb N^d$. Denote by $T:=K[y_1,\ldots,y_d]$ a standard graded polynomial ring, that is, each $y_i$ has degree $1$. We set $\tilde \Gamma_t :=\{y^{(x-h_t)/\alpha}\mid x\in\Gamma_t\}$, where $u/\alpha:=(u_{[1]}/\alpha,\ldots,u_{[d]}/\alpha)$ and $y^{u}:=y_1^{u_{[1]}}\cdot\ldots\cdot y_d^{u_{[d]}}$ for $u\in\mathbb N^d$. We get $\tilde\Gamma_t\subset T$, and hence $I_t:= \tilde\Gamma_tT$ is a monomial ideal in $T$ for all $t=1,\ldots,f$. Moreover, since $\gcd \tilde\Gamma_t=1$ we obtain $\height I_t\geq2$ (height). See \cite[Section~2]{KHSCM}. We have 
\begin{equation}\label{Kgleichung1}
K[B] \cong \bigoplus\nolimits_{t=1}^f I_t (-\deg h_t)
\end{equation}
as $\mathbb Z$-graded $T$-modules, by \cite[Proposition~2.2~(i)]{KHSCM}; where the $T$-module structure on $K[B]$ is given by the $K$-algebra homomorphism $T\hookrightarrow K[B]$, $y_i\mapsto t^{e_i}$. This implies that $\deg K[B]=f$. Applying \cite[Theorem~13.1.6]{KBSLC} to the latter ring homomorphism and to $R\twoheadrightarrow K[B]$, we get that $a(H^i_{R_+}(K[B]))=a(H^i_{T_+}(K[B]))$, since $\alpha b\in A$ for all $b\in B$. Hence the regularity of $K[B]$ is the same as an $R$-module and as a $T$-module. Then by Equation~(\ref{Kgleichung1})
\begin{equation}\label{Kregberechnung}
\reg K[B] = \max\maxk{\reg I_t +\deg h_t\mid t=1,\ldots,f};
\end{equation}
where $\reg I_t$ denotes the regularity of $I_t$ considered as a $\mathbb Z$-graded $T$-module (see also \cite[Proposition~2.2~(ii)]{KHSCM}). Note that the regularity of $K[B]$ is independent of the field $K$ for $\dim K[B]\leq5$, by \cite[Corollary~1.4]{KUBZ} and Equation~(\ref{Kregberechnung}).

\begin{bem}

One can compute this decomposition by using the \textsc{Macaulay2} \cite{KM2} package \textsc{MonomialAlgebras} \cite{KBEN}, which has been developed by Janko B\"ohm, David Eisenbud, and the author. The implemented algorithm decomposes the ring $K[Q]$ into a direct sum of monomial ideals in $K[Q']$; provided that $Q'\subseteq Q\subseteq\mathbb N^d$ are affine semigroups such that $K[Q]$ is finite over $K[Q']$. Using this decomposition we also obtain a fast algorithm computing $\reg K[Q]$ in the homogeneous case. We note that this decomposition works more general, for further details see \cite{KBENAX}.

\end{bem}

We will now recall the definition of a sequence with $*$-property:

\begin{defi}

For an element $x\in B$ we say that a sequence $\lambda=(b_1,\ldots,b_n)$ has \emph{$*$-property} if $b_1,\ldots,b_n \in \{e_1,\ldots,e_d,a_1,\ldots,a_c\}$ and $x-b_1\in B, x-b_1-b_2\in B,\ldots,x-(\sum_{j=1}^n b_j)\in B$; we say that the \emph{length} of $\La$ is $n$. Let $\lambda=(b_1,\ldots,b_n)$ be a sequence with $*$-property of $x$; we define $x(\lambda,i):=x-(\sum_{j=1}^i b_j)$ for $i=1,\ldots,n$, and $x(\lambda,0):=x$. By $\Lambda_x$ we denote the set of all sequences with $*$-property of $x$ with length $\deg x$, with the convention $\LA_0:=\emptyset$.

\end{defi}

We have $\LA_x\not=\emptyset$ for all $x\in B\setminus\{0\}$ and $\LA_x$ is finite. For elements $x,y\in \mathbb Z^d$ we define $x\geq y$ if $x_{[k]}\geq y_{[k]}$ for all $k=1,\ldots,d$.

\begin{bem}\label{Kxdegx}

Let $\lambda=(b_1,\ldots,b_n)$ be a sequence with $*$-property of $x$. For $i,j\in\mathbb N$ with $0\leq i\leq j\leq n$ we clearly have $x(\lambda,i)\geq x(\lambda,j)$. Moreover, we get $\deg x(\lambda,i)=\deg x-i$ for $i=0,\ldots,n$. Hence in case that $\La\in\LA_x$ it follows that $x(\La,\deg x)=0$. 

\end{bem}

The definition of a sequence with $*$-property is motivated to control the degree of $K[B]$ in a certain way:

\begin{lem}[{\cite[Lemma~2.4]{KMNSNAX}}]\label{K*prop}

Let $x\in B_A\setminus\{0\}$ and $\lambda=(b_1,\ldots,b_n)$ be a sequence with $*$-property of $x$. Then

\begin{enumerate}

\item[$1)$] $x(\lambda,i)\in B_A$ for all $i=0,\ldots,n$.

\item[$2)$] $x(\lambda,i)\not\sim x(\lambda,j)$ for all $i,j\in\mathbb N$ with $0\leq i<j\leq n$.

\end{enumerate}

\end{lem}

Thus, we get $\deg x\leq \deg K[B]-1$ for every $x\in B_A$.

\begin{bem}\label{Ka1bisc}

By construction we have $\{0,a_1,\ldots,a_c\}\subseteq B_A$. Let $x\in\{0,a_1,\ldots,a_c\}$ and $y\in B_A$ with $x\not=y$, and suppose that $x\sim y$. Since $0\leq x_{[i]}<\alpha$ for all $i=1,\ldots,d$ we obtain $y\geq x$, and therefore $y-x\in A\setminus\{0\}$ which contradicts $y\in B_A$. Thus, we have $x\not\sim y$. Hence one could assume that $\Gamma_1=\{0\},\Gamma_2=\{a_1\},\ldots,\Gamma_{c+1}=\{a_c\}$.

\end{bem}

Combining Lemma~\ref{K*prop} and Remark~\ref{Ka1bisc} one can show that (see also \cite[Proposition~2.7]{KMNSNAX}):

\begin{prop}[{\cite[Theorem~1.1]{KHSCM}}]\label{Kegred}

We have $\deg x\leq \deg K[B] - \codim K[B]$ for all $x\in B_A$.

\end{prop}

The next definition will be useful to count the number of equivalence classes which arises from two sequences with $*$-property.

\begin{defi}

Let $x,y\in B\setminus\{0\}$, $\lambda\in\Lambda_x$, and $\nu\in\Lambda_y$. We define

\begin{enumerate}

\item[$1)$] $\Delta(\lambda,\nu):=\{(i,j)\in\mathbb N^2\mid 0\leq i\leq\deg x, 0\leq j\leq \deg y, x(\lambda,i)\sim y(\nu,j)\}$,

\item[$2)$] $\delta(\lambda,\nu):=\#\Delta(\lambda,\nu)-2$, and

\item[$3)$] $\delta(x,y):=\min\limits_{\lambda'\in\Lambda_x,\,\nu'\in\Lambda_y} \delta(\lambda',\nu')$.

\end{enumerate}

\end{defi}

For example, $\Delta((e_1),(e_2))=\{(0,0),(0,1),(1,0),(1,1)\}$, where $(e_1)\in\LA_{e_1}$ and $(e_2)\in\LA_{e_2}$, hence $\delta((e_1),(e_2))=2$. Moreover, since $\#\LA_{e_1}=\#\LA_{e_2}=1$ we obtain $\delta(e_1,e_2)=2$.

\begin{bem}\label{Kprojektion}

Let $x,y\in B\setminus\{0\}$ with $x\sim y$, $\lambda\in\Lambda_x$, and $\nu\in\Lambda_y$. We always have $(0,0),(\deg x,\deg y)\in\Delta(\La,\n)$, since $x(\lambda,0)\sim y(\nu,0)$ and $x(\lambda,\deg x)\sim y(\nu,\deg y)$. Hence $\delta(\La,\n)\geq 0$ and $\delta(x,y)\geq0$. Moreover, if $x,y\in B_A$ and $(i,j)\in \Delta(\La,\n)$, then $(i,k)\notin\Delta(\La,\n)$ for all $k\in\{0,\ldots,\deg y\}\setminus\{j\}$ by Lemma~\ref{K*prop}~$2)$, since otherwise $y(\nu,j)\sim y(\nu,k)$ for $j\not=k$. This argument shows that $\#\Delta(\La,\n)\leq\min\maxk{\deg x,\deg y}+1$ in case that $x,y\in B_A$.

\end{bem}

Consider the set $L=\{x(\La,0),\ldots,x(\La,\deg x)\}\cup\{y(\nu,0),\ldots,y(\nu,\deg y)\}$, where $\La\in\LA_x$ and $\nu\in\LA_y$. Since $\deg K[B]$ is equal to the number of equivalence classes on $B$ and on $B_A$ it would be nice to know the number of equivalence classes on $L$. The latter number is equal to $\deg x+\deg y-\delta(\La,\nu)$ if $x,y\in B_A\setminus\{0\}$, by Lemma~\ref{K*prop}~$2)$. Thus, it is of interest to have a good bound for $\delta(\La,\nu)$ and for $\delta(x,y)$; to state our upper bound conjecture we need one more definition.

\begin{defi}

Let $x,y\in B$. We define $h(x,y)$ by:
$$
h(x,y):=(\min\{x_{[1]},y_{[1]}\}, \min\{x_{[2]},y_{[2]}\},\ldots,\min\{x_{[d]},y_{[d]}\}).
$$

\end{defi}

In case $x\sim y$ we get that $h(x,y)\sim x,y$, and hence $h(x,y)\in G(B)\cap\mathbb N^d$. We conjecture

\begin{con}\label{Kkomisch}

Let $x,y\in B_A\setminus\{0\}$ with $x\sim y$. Then 
$$
\delta(x,y)\leq \deg h(x,y)-1.
$$

\end{con}

Conjecture~\ref{Kkomisch} holds in case that $x=y$, since $\delta(x,x)\leq\delta(\La,\nu)\leq \deg x-1=\deg h(x,x)-1$ for all $\La,\nu\in\LA_x$, by Remark~\ref{Kprojektion}. The inequality in the conjecture clearly does not hold for equivalent elements of $B$, for example, $\delta(e_1,e_2)=2>-1=\deg h(e_1,e_2)-1$. Moreover, one could easily construct a stronger version of Conjecture~\ref{Kkomisch}, namely, for all $x,y\in B_A\setminus\{0\}$ with $x\sim y$, and for each pair $\La\in\LA_x$ and $\nu\in\LA_y$ we have $\delta(\La,\nu)\leq \deg h(x,y)-1$.

\begin{bsp}\label{Kbspcross1}

Consider the monoid $B=\erz{(30,0),(0,30),(3,27),(23,7)}$, $x=(27,243)$, and $y=(207,63)$. We get that $x,y\in B_A$ with $x\sim y$, since $x-y=(-180,180)\in G(A)=30\mathbb Z^2$. Clearly, $\Lambda_x=\{((3,27),\ldots,(3,27))\}=\{\lambda\}$ and $\Lambda_y=\{((23,7),\ldots,(23,7))\}=\{\nu\}$. We have $\delta(x,y)=2$, since $\Delta(\La,\n)=\{(0,0), (3,3), (6,6), (9,9)\}$ and $\#\Lambda_x=\#\Lambda_y=1$. Moreover, $\deg h(x,y)=\deg{(27,63)}=3$ and therefore Conjecture~\ref{Kkomisch} holds and is sharp. 

\end{bsp}

\begin{nota}

Let $x\in B\setminus\{0\}$.  It is often useful to illustrate a sequence with $*$-property $\lambda\in\Lambda_x$ as a graph, where the set of vertices is a subset of $\{x(\lambda,i)\mid i\in\{0,\ldots,\deg x\}\}$. Let $x(\La,i)$ and $x(\La,j)$ be vertices; there will be an edge between $x(\La,i)$ and $x(\La,j)$ if $j>i$ and there is no vertex $x(\La,k)$ with $j>k>i$. Moreover, $x$ and $0$ will always be vertices, keep in mind that $x(\La,0)=x$ and $x(\La,\deg x)=0$. So Example~\ref{Kbspcross1} can be illustrated by the graph 
$$
\begin{xy}
  \xymatrix{
x(\La,0)=x \ar@{-}[rr] &  & x(\lambda,3)\ar@{-}[rr]& & x(\lambda,6)\ar@{-}[rr]  &  & x(\lambda,9) = 0,
  }
\end{xy}
$$
and by the graph 
$$
\begin{xy}
  \xymatrix{
y(\nu,0)=y \ar@{-}[rr] & & y(\nu,3)\ar@{-}[rr]& & y(\nu,6)  \ar@{-}[rr]  & & y(\nu,9) = 0.
  }
\end{xy}
$$
To get a better understanding and to avoid extensive writing we illustrate these situations by a picture:
$$
\begin{xy}
  \xymatrix{
   x \ar@{-}[rr] \ar@{~}[d] &  & x(\lambda,3)\ar@{-}[rr]\ar@{~}[d] & & x(\lambda,6)\ar@{-}[rr] \ar@{~}[d] &  & 0 \ar@{~}[d]\\
y \ar@{-}[rr] & & y(\nu,3)\ar@{-}[rr]& & y(\nu,6)  \ar@{-}[rr]  & & \,0,
  }
\end{xy}
$$
where the sidled lines denote equivalent elements.  \textbf{Sidled lines always denote equivalent elements, though equivalent elements may not be illustrated in such a picture.}

\end{nota}

\begin{defi}

Let $x,y\in B\setminus \{0\}$, $\lambda\in\Lambda_x$, and $\nu\in\Lambda_y$.

\begin{enumerate}

\item[$1)$] Let $(i,j),(i',j')\in\Delta(\La,\n)$. We define a partial order $\leq$ on $\Delta(\La,\n)$ by $(i,j)\leq (i',j')$ if $i\leq i'$ and $j\leq j'$.

\item[$2)$] We say that $\lambda$ and $\nu$ are \emph{crossless} if $(\Delta(\La,\n),\leq)$ is a totally ordered set, meaning for all $(i,j),(i',j')\in\Delta(\La,\n)$ we have $(i,j)\leq(i',j')$ or $(i,j)\geq (i',j')$.

\item[$3)$] We say that $x$ and $y$ are \emph{crossless} if there exist sequences with $*$-property $\lambda'\in\Lambda_x$ and $\nu'\in\Lambda_y$ which are crossless.

\end{enumerate}

\end{defi}

Thus, in case that two sequences with $*$-property $\La\in\LA_x$ and $\nu\in\LA_y$ are not crossless there exist $i,j,k,l\in\mathbb N$ with $i<j\leq\deg x$ and $l<k\leq\deg y$ such that $x(\La,i)\sim y(\nu,k)$ and $x(\La,j)\sim y(\nu,l)$. We note that $\La$ and $\nu$ in Example~\ref{Kbspcross1} are crossless. Moreover, for $x\in B_A\setminus\{0\}$ we get that $\La$ and $\La$ are crossless for all $\La\in \LA_x$, by Lemma~\ref{K*prop}~$2)$. This property is of interest, since Conjecture~\ref{Kkomisch} as well as the stronger version hold for crossless elements, see Proposition~\ref{Khungleichung} and Corollary~\ref{Kcrosslesskor}. Our original thought was that equivalent elements in $B_A$ should be crossless, however, this is not the case:

\begin{bsp}\label{Kkomischesbeispiel}

Consider the monoid $B=\erz{(79,0),(0,79),(77,2),(34,45)}$, $x=(1232,32)$, and $y=(442,585)$. We have $x,y\in B_A$ with $x\sim y$. We get $\Lambda_x=\{((77,2),\ldots,(77,2))\}=\{\lambda\}$, and $\Lambda_y=\{((34,45),\ldots,(34,45))\}=\{\nu\}$. Moreover, $\Delta(\La,\n)=\{(0,0),(5,9),(11,4),(16,13)\}$ and therefore $x$ and $y$ are not crossless. This situation can be illustrated by:
$$
\begin{xy}
  \xymatrix{
x \ar@{-}[rr] \ar@{~}[d] &  & x(\lambda,5)\ar@{-}[rr]\ar@{~}[drr] & & x(\lambda,11)\ar@{-}[rr] &  & 0 \ar@{~}[d]\\
y \ar@{-}[rr] & & y(\nu,4)\ar@{-}[rr]\ar@{~}[urr] & & y(\nu,9)  \ar@{-}[rr]  & & \,0.}
\end{xy}
$$
We get $\delta(\La,\n)=\delta(x,y)=2$ and $\deg h(x,y)=\deg {(442,32)}=6$, hence Conjecture~\ref{Kkomisch} holds.

\end{bsp}

\begin{bem}\label{Kpermutation}

We note that $(b_{\sigma(1)},\ldots,b_{\sigma(n)})$ is a sequence with $*$-property of $x$ for every sequence with $*$-property $\La=(b_1,\ldots,b_{n})$ of $x$ and for every permutation $\sigma$ of $\{1,\ldots,n\}$, since $x=x(\La,n) + \sum_{j=1}^nb_j$. This leads to:

\end{bem}

\begin{defi}

Let $\lambda=(b_1,b_2,\ldots,b_{n})$ be a sequence with $*$-property of $x$. We define $\La^*:=(b_{n},b_{n-1},\ldots,b_1)$ as the trivial permutation of $\lambda$.

\end{defi}

Thus, by construction $\La\in\LA_x$ if and only if $\La^*\in\LA_x$. Moreover, for a sequence with $*$-property $\La=(b_1,\ldots,b_{\deg x})\in\LA_x$ and $i\in\{0,\ldots,\deg x\}$ we have the following symmetry:
\begin{equation}\label{K*prop2}
x-x(\La,i)=x-(x-\sum_{j=1}^ib_j)=\sum_{j=1}^ib_j=x-\sum_{j=1}^{\deg x-i}b_{\deg x+1-j}=x(\lambda^*,\deg x-i).
\end{equation}

\begin{bem}\label{Kaussehen1}

Let $x\in B\setminus\{0\}, \lambda=(b_1,\ldots,b_{\deg x})\in\Lambda_x$, and $i\in\{1,\ldots,\deg x-1\}$. We have $(b_1,\ldots,b_i)\in\Lambda_{x(\lambda^*,\deg x-i)}$, since $x(\lambda^*,\deg x-i)=\sum_{j=1}^ib_j$; moreover, we obtain $(b_{i+1},\ldots,b_{\deg x})\in\Lambda_{x(\lambda,i)}$, since $x(\lambda,i)=\sum_{j=1}^{\deg x-i} b_{i+j}$.

\end{bem}

For sets $X,Y\subseteq \mathbb N^d$ we define the set $X+Y:=\{x+y\mid x\in X,y\in Y\}\subseteq\mathbb N^d$. We will now establish two basic combinatorial Lemmas which are essential to study the crossless property:

\begin{lem}\label{Klemeqcross}

Let $x,y\in B\setminus\{0\}$ with $x\sim y$. Moreover, let $\lambda=(b_1,\ldots,b_{\deg x})\in\Lambda_x$ and $\nu=(g_1,\ldots,g_{\deg y})\in\Lambda_y$. Assume that there is some $i\in\{1,\ldots,\deg x-1\}$ and some $k\in\{1,\ldots,\deg y-1\}$ such that $x(\La,i)\sim y(\nu,k)$, that is,
$$
\begin{xy}
  \xymatrix{
x \ar@{-}[rr]^{} \ar@{~}[d] & & x(\lambda,i)\ar@{-}[rr]^{}\ar@{~}[d] &  & 0 \ar@{~}[d]\\
y \ar@{-}[rr]^{} & & y(\nu,k)  \ar@{-}[rr]^{}  & & \,0.
  }
\end{xy}
$$
Set $x':=x(\lambda^*,\deg x-i)$, $x'':=x(\lambda,i)$, $y':=y(\nu^*,\deg y-k)$, and $y'':=y(\nu,k)$. Moreover, set $\lambda':=(b_1,\ldots,b_i)\in\Lambda_{x'}$, $\lambda'':=(b_{i+1},\ldots,b_{\deg x})\in\Lambda_{x''}$, $\nu':=(g_1,\ldots,g_k)\in\Lambda_{y'}$, and $\nu'':=(g_{k+1},\ldots,g_{\deg y})\in\Lambda_{y''}$ (see Remark~\ref{Kaussehen1}). We have

\begin{enumerate}

\item[$1)$] $x(\lambda^*,\deg x-i)\sim y(\nu^*,\deg y-k).$ 

\item[$2)$] $\Delta(\La',\n')=\{(m,n)\in\Delta(\La,\n)\mid (m,n)\leq (i,k)\}$.

\item[$3)$] $\{(i,k)\} + \Delta(\La'',\n'')=\{(m,n)\in\Delta(\La,\n)\mid (m,n)\geq (i,k)\}$.

\item[$4)$] If $\lambda$ and $\nu$ are crossless, then $\lambda'$ and $\nu'$ are crossless.

\item[$5)$] If $\lambda$ and $\nu$ are crossless, then $\La''$ and $\n''$ are crossless.

\item[$6)$] $\delta(\lambda',\nu') + \delta(\lambda'',\nu'') \leq \delta(\lambda,\nu)-1.$
Equality holds, if $\lambda$ and $\nu$ are crossless.

\end{enumerate}

\end{lem}
\begin{proof}

$1)$ We have $x-x(\La,i)=x(\lambda^*,\deg x-i)$ and $y-y(\nu,k)=y(\nu^*,\deg y-k)$ by Equation~(\ref{K*prop2}). The assertion follows from $x-y, y(\nu,k)-x(\lambda,i)\in G(A)$.\\
$2)$ Note that $\Delta(\lambda',\nu')\subseteq\{0,\ldots,i\}\times\{0,\ldots,k\}$. Let $m,n\in\mathbb N$ with $m\leq i$ and $n\leq k$. We get $x(\La,m)-x'(\La',m)=x(\La,i)$ as well as $y(\n,n)-y'(\n',n)=y(\n,k)$. Hence
$$
x(\La,m)-y(\n,n) + y'(\n',n)-x'(\La',m)\in G(A),
$$
which proves $2)$.\\
$3)$ Note that $\Delta(\lambda'',\nu'')\subseteq\{0,\ldots,\deg x-i\}\times\{0,\ldots,\deg y-k\}$. Let $m,n\in\mathbb N$ with $m\leq\deg x-i$ and $n\leq\deg y-k$. The assertion follows from
$$
x''(\La'',m)=x(\lambda,i+m)\quad\mbox{and}\quad y''(\n'',n)=y(\n,k+n).
$$
$4), 5)$ This follows from $2)$ and $3)$.\\
$6)$ Since $(i,k)\in\Delta(\La',\n'),(0,0)\in\Delta(\lambda'',\nu'')$, and  $\Delta(\La',\n')\subseteq\{0,\ldots,i\}\times\{0,\ldots,k\}$, we have 
\begin{equation}\label{Kgleichh11}
\#\left(\Delta(\lambda',\n')\cap(\{(i,k)\} + \Delta(\La'',\n''))\right)=1.
\end{equation}
Hence
\begin{equation}\label{Keq1}
\#\Delta(\lambda',\n')+\#\Delta(\lambda'',\nu'')-1\stackrel{\mbox{\scriptsize{(\ref{Kgleichh11})}}}{=}\#\left(\Delta(\lambda',\n')\cup(\{(i,k)\} + \Delta(\La'',\n''))\right)\stackrel{\mbox{\scriptsize{2),3)}}}{\leq}\#\Delta(\La,\n),
\end{equation}
and therefore
\begin{equation}\label{Keq3343}
\delta(\lambda',\nu') + \delta(\lambda'',\nu'')=\#\Delta(\lambda',\n')+\#\Delta(\lambda'',\nu'')-1-3\stackrel{\mbox{\scriptsize{(\ref{Keq1})}}}{\leq}\#\Delta(\La,\n)-2-1=\delta(\lambda,\nu)-1.
\end{equation}
If $\La$ and $\n$ are crossless we have equality in Equation~(\ref{Keq1}), by assertion $2)$ and $3)$. Hence we also have equality in Equation~(\ref{Keq3343}).
\end{proof}

\begin{bem}\label{KmingensK}

Consider two elements $x,y\in B_A$ with $x\sim y$ and $x\not=y$. Suppose that $x\leq y$ or $x\geq y$, hence $y-x\in A\setminus\{0\}$ or $x-y\in A\setminus\{0\}$ which contradicts $y\in B_A$ or $x\in B_A$. Thus, there exists $i,j\in\{1,\ldots,d\}$ such that $x_{[i]}>y_{[i]}$ and $x_{[j]}<y_{[j]}$. This argument also shows that $\tilde\Gamma_t$ is a minimal generating set of $I_t$ for all $t=1,\ldots,f$.

\end{bem}

\begin{lem}\label{Kcrosslem}

Let $x,y\in B_A\setminus\{0\}$ with $x\sim y$, $\lambda\in\Lambda_x$, and $\nu\in\Lambda_y$. Moreover, let $\La$ and $\n$ be not crossless, that is, $x(\La,i)\sim y(\nu,k)$ and $x(\La,j)\sim y(\nu,l)$ for some $i,j,l,k\in\mathbb N$ with $i<j\leq\deg x$ and $l<k\leq \deg y$. We have:

\begin{enumerate}

\item[$1)$] $i\geq2$, $l\geq2$, $j\leq\deg x-2$, and $k\leq\deg y-2$. This can be illustrated by:
$$
\begin{xy}
  \xymatrix{
x \ar@{-}[rr]^{} \ar@{~}[d] &  & x(\lambda,i)\ar@{-}[rr]^{}\ar@{~}[drr] & & x(\lambda,j)\ar@{-}[rr]^{} &  & 0 \ar@{~}[d]\\
y \ar@{-}[rr]^{} & & y(\nu,l)\ar@{-}[rr]^{}\ar@{~}[urr] & & y(\nu,k)  \ar@{-}[rr]^{}  & & \,0.
  }
\end{xy}
$$

\item[$2)$] $\lambda^*$ and $\nu^*$ are not crossless, in particular:
$$
\begin{xy}
  \xymatrix{
x \ar@{-}[rr]^{} \ar@{~}[d] &  & x(\lambda^*,\deg x-j)\ar@{-}[rr]^{}\ar@{~}[drr] & & x(\lambda^*,\deg x-i)\ar@{-}[rr]^{} &  & 0 \ar@{~}[d]\\
y \ar@{-}[rr]^{} & & y(\nu^*,\deg y-k)\ar@{-}[rr]^{}\ar@{~}[urr] & & y(\nu^*,\deg y-l)  \ar@{-}[rr]^{}  & & \,0.
  }
\end{xy}
$$

\item[$3)$] $x(\La,i)\not=y(\n,k)$ and $x(\La,j)\not=y(\n,l)$.

\item[$4)$] $y(\n,k)_{[n]}>x(\La,i)_{[n]}$ and $x(\La,j)_{[m]}>y(\n,l)_{[m]}$ for some $n,m\in\{1,\ldots,d\}$ with $n\not=m$.

\item[$5)$] $y(\n,k)_{[n']}<x(\La,i)_{[n']}$ and $x(\La,j)_{[m']}<y(\n,l)_{[m']}$ for some $n',m'\in\{1,\ldots,d\}$.

\end{enumerate}

\end{lem}

\begin{proof}

$1),2)$ We have $i,l\not=0$, $j\not=\deg x$, as well as $k\not=\deg y$ by Lemma~\ref{K*prop}~$2)$. By Lemma~\ref{Klemeqcross}~$1)$ we get $x(\lambda^*,\deg x-i)\sim y(\nu^*,\deg y-k)$ and $x(\lambda^*,\deg x-j)\sim y(\nu^*,\deg y-l)$. Since $\deg x-i>\deg x-j$ and $\deg y-k<\deg y-l$, it follows that $\lambda^*$ and $\n^*$ are not crossless which proves assertion $2)$. Suppose that $j=\deg x-1$, that is, $\deg x(\La,j)=1$. We have $y(\nu,l)\in B_A$ by Lemma~\ref{K*prop}~$1)$ and $y(\nu,l)\sim x(\La,j)$, hence $y(\nu,l)= x(\La,j)$ by Remark~\ref{Ka1bisc}, which contradicts $l<k<\deg y$. Assertion $1)$ now follows from symmetry and $2)$.\\
$3)$ By symmetry we only need to show that $x(\lambda,i)\not= y(\nu,k)$. Let $\lambda=(b_1,\ldots,b_{\deg x})$ and $\nu=(g_1,\ldots,g_{\deg y})$ and set $\nu':=(g_1,\ldots,g_k,b_{i+1},\ldots,b_{\deg x})$. Suppose to the contrary that $x(\lambda,i)=y(\nu,k)$. We get $\nu'\in\LA_y$, moreover, $y(\nu',k+j-i)=x(\lambda,j)\sim y(\nu,l)=y(\nu',l)$, which contradicts Lemma~\ref{K*prop}~$2)$, since $k+j-i>l$.\\ 
$4),5)$ By Lemma~\ref{K*prop}~$1)$ we have $x(\La,i),y(\n,k)\in B_A$. Moreover, $x(\La,i)\not=y(\n,k)$ by assertion~$3)$ and $x(\La,i)\sim y(\nu,k)$ by assumption. Hence $y(\n,k)_{[n]}>x(\La,i)_{[n]}$ as well as $y(\n,k)_{[n']}<x(\La,i)_{[n']}$ for some $n,n'\in\{1,\ldots,d\}$ by Remark~\ref{KmingensK}. Analogously, $x(\La,j)_{[m]}>y(\n,l)_{[m]}$ and $x(\La,j)_{[m']}<y(\n,l)_{[m']}$ for some $m,m'\in\{1,\ldots,d\}$. Suppose to the contrary that $m=n$; then $x(\La,j)_{[m]}>y(\n,l)_{[m]}\geq y(\n,k)_{[m]}>x(\La,i)_{[m]}\geq x(\La,j)_{[m]}$, which is a contradiction.
\end{proof}

\begin{lem}\label{Kcrossgeq}

Consider the same situation as in Lemma~\ref{Kcrosslem}. Let $n,m\in\{1,\ldots,d\}$ such that $y(\n,k)_{[n]}>x(\La,i)_{[n]}$ and $x(\La,j)_{[m]}>y(\n,l)_{[m]}$. Then

\begin{enumerate}

\item[$1)$] $y(\n,l)_{[n]}>x(\La,j)_{[n]}$.

\item[$2)$] $x(\La,i)_{[m]}>y(\n,k)_{[m]}$.

\end{enumerate}

\end{lem}
\begin{proof}

$1)$ We have $y(\n,l)_{[n]}\geq y(\n,k)_{[n]}>x(\La,i)_{[n]}\geq x(\La,j)_{[n]}$.\\
$2)$ We have $x(\La,i)_{[m]}\geq x(\La,j)_{[m]}>y(\n,l)_{[m]}\geq y(\n,k)_{[m]}$.
\end{proof}

Considering again Example~\ref{Kkomischesbeispiel} one can show that $(521,190)\in B_A$, moreover, it is equivalent to $(1232,32)$ and to $(442,585)$. This holds in general:

\begin{prop}\label{Keinsnoch}

Let $x,y\in\Gamma_t\subseteq{B_A}\setminus\{0\}$ for some $t\in\{1,\ldots,f\}$, $\lambda\in\Lambda_x$, and $\nu\in\Lambda_y$. If $\lambda$ and $\nu$ are not crossless, then there is some $z\in\Gamma_t$ with $z\not=x$ and $z\not=y$.

\end{prop}

\begin{proof}

By Lemma~\ref{Kcrosslem} we have
$$
\begin{xy}
  \xymatrix{
x \ar@{-}[rr]^{} \ar@{~}[d] &  & x(\lambda,i)\ar@{-}[rr]^{}\ar@{~}[drr] & & x(\lambda,j)\ar@{-}[rr]^{} &  & 0 \ar@{~}[d]\\
y \ar@{-}[rr]^{} & & y(\nu,l)\ar@{-}[rr]^{}\ar@{~}[urr] & & y(\nu,k)  \ar@{-}[rr]^{}  & & \,0,
  }
\end{xy}
$$
for some $i,j,l,k\in\mathbb N$ with $0<i<j< \deg x$ and $0<l<k< \deg y$. We set 
$$
z':=x(\lambda,j)+y(\nu^*,\deg y-l).
$$ 
By construction we have $z'\in B$. By Lemma~\ref{Kcrosslem}~$5)$ we get
$$
x(\lambda,j)_{[h]}<y(\nu,l)_{[h]}
$$
for some $h\in\{1,\ldots,d\}$. Hence $z'_{[h]}<y_{[h]}$ by Equation~(\ref{K*prop2}). Applying Lemma~\ref{Kcrosslem}~$5)$ to the second assertion in the same lemma we obtain
 $$
 y(\nu^*,\deg y-l)_{[g]}< x(\lambda^*,\deg x-j)_{[g]}
 $$
for some $g\in\{1,\ldots,d\}$. Hence $z'_{[g]}<x_{[g]}$ again by Equation~(\ref{K*prop2}). Consider an element $z:=z'-\sum_{u=1}^dn_ue_u\in B$ for some $n_u\in\mathbb N$ such that $\sum_{u=1}^{d} n_u$ is maximal. We obtain $z\in B_A$, in particular $z\leq z'$ and therefore $z\not=x,y$. Moreover, by construction $z\sim z'$. We get 
 $$
 z'-x=x(\lambda,j)+y(\nu^*,\deg y-l)-x\stackrel{\mbox{\scriptsize{(\ref{K*prop2})}}}{=}y(\nu^*,\deg y-l)-x(\lambda^*,\deg x-j)\stackrel{\mbox{\scriptsize{\ref{Klemeqcross}}}}{\in} G(A).
 $$
Hence $z'\sim x$, that is, $z\in\Gamma_t$.
\end{proof}

\begin{fol}\label{Kcross2el}

Let $\#\Gamma_t=2$ for some $t\in\{1,\ldots,f\}$, say $\Gamma_t=\{x,y\}$, $\La\in\LA_x$, and $\n\in\LA_y$. Then $\La$ and $\n$ are crossless, in particular $x$ and $y$ are crossless.

\end{fol}

\begin{proof}

Suppose that $\La$ and $\n$ are not crossless. By Proposition~\ref{Keinsnoch} we get some element $z\in\Gamma_t$ with $z\not=x,y$, which contradicts $\#\Gamma_t=2$. Hence $x$ and $y$ are crossless as well.
\end{proof}

\begin{bem}\label{Kungleichungh}

Let $x',x'',y',y''\in B$. Then 
$$
h(x',y')+h(x'',y'')\leq h(x'+x'',y'+y'').
$$ 
This follows immediately from the definition, however, here is a
\begin{proof}
Let $i\in\{1,\ldots,d\}$. We have
$$
2\min\maxk{(x'+x'')_{[i]},(y'+y'')_{[i]}}=x'_{[i]}+y'_{[i]}+x''_{[i]}+y''_{[i]}-|x'_{[i]}-y'_{[i]}+x''_{[i]}-y''_{[i]}|\hspace{1.35cm}
$$
$$
\geq x'_{[i]}+y'_{[i]}-|x'_{[i]}-y'_{[i]}|+x''_{[i]}+y''_{[i]}-|x''_{[i]}-y''_{[i]}|=2\min\maxk{x'_{[i]},y'_{[i]}}+2\min\maxk{x''_{[i]},y''_{[i]}}.
$$
and therefore $h(x',y')_{[i]}+h(x'',y'')_{[i]}\leq h(x'+x'',y'+y'')_{[i]}$ and we are done.
\end{proof}
\end{bem}

As a consequence of the next proposition Conjecture~\ref{Kkomisch} holds for crossless elements. Moreover, Proposition~\ref{Khungleichung} together with Corollary~\ref{Kcross2el} leads to a proof of the Eisenbud-Goto conjecture for the at most two element case in Section~\ref{Kcase2}.

\begin{prop}\label{Khungleichung}

Let $x,y\in B_A\setminus\{0\}$ with $x\sim y$, $\lambda\in\LA_x$, and $\n\in\LA_y$. If $\lambda$ and $\n$ are crossless, then 
$$
\delta(\La,\n)\leq \deg h(x,y)-1.
$$

\end{prop}
\begin{proof}

We show this by induction on $\delta(\La,\n)\in\mathbb N$. Let $\delta(\La,\n)=0$, that is, we need to show that $\deg h(x,y)\geq 1$. Suppose that $\deg h(x,y)<1$, that is, $h(x,y)=0$, since $h(x,y)\in G(B)\cap\mathbb N^d$. Hence $x,y\sim0$, since $h(x,y)\sim x,y$, which contradicts $x,y\in B_A\setminus\{0\}$. In case that $\delta(\lambda,\n)>0$ we can fix an $i\in\{1,\ldots,\deg x-1\}$ such that $x(\La,i)\sim y(\n,k)$ for some $k\in\{1,\ldots\deg y-1\}$ by Lemma~\ref{K*prop}~$2)$. Using the notation of Lemma~\ref{Klemeqcross} as well as assertion $4)$ and $5)$ we get that $\La'\in\LA_{x'}$ and $\n'\in \LA_{y'}$ are crossless, and also that $\La''\in\LA_{x''}$ and $\n''\in \LA_{y''}$ are crossless, since $\La$ and $\n$ are crossless. We have $x',x'',y',y''\in B_A\setminus\{0\}$ by Lemma~\ref{K*prop}~$1)$, moreover, $x''\sim y''$ and $x'\sim y'$ by Lemma~\ref{Klemeqcross}~$1)$. Note that $x=x'+x''$ and $y=y'+y''$ by Equation~(\ref{K*prop2}). Hence by induction
$$
\delta(\La,\n)\stackrel{\mbox{\scriptsize{\ref{Klemeqcross}}}}{=}\delta(\La',\n')+\delta(\La'',\n'')+1\leq \deg h(x',y')+\deg h(x'',y'')-1\stackrel{\mbox{\scriptsize{\ref{Kungleichungh}}}}{\leq} \deg h(x,y)-1.
$$
\end{proof}

\begin{fol}\label{Kcrosslesskor}

Let $x,y\in B_A\setminus\{0\}$ with $x\sim y$. If $x$ and $y$ are crossless, then 
$$
\delta(x,y)\leq \deg h(x,y)-1.
$$

\end{fol}
\begin{proof}

Since $x$ and $y$ are crossless, there are some sequences $\lambda\in\LA_x$ and $\n\in\LA_y$ which are crossless. The assertion follows from $\delta(x,y)\leq \delta(\lambda,\n)$ and Proposition~\ref{Khungleichung}.
\end{proof}

In Proposition~\ref{Kkomischclose} we will prove Conjecture~\ref{Kkomisch} for certain elements in the dimension $2$ case. In the proof of this proposition we can assume that these elements are not crossless by the above corollary; the next definition is motivated to get a better handling of this situation.

\begin{defi}

Let $x,y\in B\setminus\{0\}$. By a \emph{cross} we mean a tuple $(\La,\n,i,j,l,k)\in\LA_x\times\LA_y\times\mathbb N^4$ with $i<j\leq\deg x$ and $l<k\leq\deg y$ such that $x(\La,i)\sim y(\n,k)$ and $x(\La,j)\sim y(\n,l)$. In this case we will say that $\La$ and $\nu$ have a cross. Moreover, the \emph{height} of a cross $(\La,\n,i,j,l,k)$ is defined to be $(j-i,k-l)\in\mathbb N^2$.

\end{defi}

Thus, the sequences $\La$ and $\nu$ are not crossless if and only if they have a cross. To end this section consider the following situation: let $x,y\in B\setminus\{0\}$ with $x\sim y$, $\La\in\LA_x$, and $\nu\in\LA_y$, moreover, let $(\La,\n,i,j,l,k)$ and $(\La,\n,i',j',l',k')$ be crosses with $j\leq i'$ and $k\leq l'$. This can be illustrated by
$$
\begin{xy}
  \xymatrix{
x \ar@{-}[r]^{} \ar@{~}[d] & x(\lambda,i)\ar@{-}[r]^{}\ar@{~}[dr] & x(\lambda,j)\ar@{-}[r]^{} & x(\lambda,i')\ar@{-}[r]^{}\ar@{~}[dr] & x(\lambda,j')\ar@{-}[r]^{} & 0 \ar@{~}[d]\\
y \ar@{-}[r]^{} & y(\nu,l)\ar@{-}[r]^{}\ar@{~}[ur] & y(\nu,k)  \ar@{-}[r]^{} & y(\nu,l')\ar@{-}[r]^{}\ar@{~}[ur] & y(\nu,k')  \ar@{-}[r]^{} & \,0;
  }
\end{xy}
$$
assume in the picture $0<i<j<i'<j'<\deg x$ and $0<l<k<l'<k'<\deg y$. In Proposition~\ref{Kkomischclose} it will be very important to prevent such a case if one cross is fixed. We can avoid this, for example, by assuming that $j-i$ (or $j'-i'$) is maximal among all crosses:

\begin{prop}\label{Kgluingcrosses}

Let $x,y\in B\setminus\{0\}$ with $x\sim y$, $\La\in\LA_x$, and $\n\in\LA_y$. Moreover, let $(\La,\n,i,j,l,k)$ and $(\La,\n,i',j',l',k')$ be crosses such that $j\leq i'$ and $k\leq l'$. Then there is a cross $(\La',\n',i'-j,j'-i,l'-k,k'-l)\in\LA_x\times\LA_y\times \mathbb N^4$ of height $(j-i+j'-i',k-l+k'-l')$.

\end{prop}
\begin{proof}

Let $\La=(b_1,\ldots,b_{\deg x})$ and $\n=(g_1,\ldots,g_{\deg y})$. Set 
$$
\La':=(b_{j+1},\ldots,b_{j'},b_{i+1},\ldots,b_{j},b_{1},\ldots,b_{i},b_{j'+1},\ldots,b_{\deg x})
$$
and
$$
\n':=(g_{k+1},\ldots,g_{k'},g_{l+1},\ldots,g_{k},g_{1},\ldots,g_{l},g_{k'+1},\ldots,g_{\deg y}).
$$
By Remark~\ref{Kpermutation} we obtain $\La'\in\LA_x$ and $\n'\in\LA_y$. We have $i'-j<j'-i$ and $k'-l>l'-k$. We claim that $x(\La',i'-j)\sim y(\n',k'-l)$ and $x(\La',j'-i)\sim y(\n',l'-k)$, and therefore $(\La',\n',i'-j,j'-i,l'-k,k'-l)$ is a cross of height $(j-i+j'-i', k-l+k'-l')$. To verify the claim, note that
$$
x(\La',i'-j)=x-\sum_{t=1}^{i'-j}b_{j+t}=x-(x(\La,j)-x(\La,i'))\sim y-(y(\n,l)-y(\n,k'))
$$
$$
=y-\sum_{t=1}^{k'-l}g_{l+t}=y-\sum_{t=1}^{k'-k}g_{k+t}-\sum_{u=1}^{k-l}g_{l+u}=y(\n',k'-l),\hspace{3.55cm}
$$
and
$$
y(\n',l'-k)=y-\sum_{t=1}^{l'-k}g_{k+t}=y-(y(\n,k)-y(\n,l'))\sim x-(x(\La,i)-x(\La,j'))
$$
$$
=x-\sum_{t=1}^{j'-i}b_{i+t}=x-\sum_{t=1}^{j'-j}b_{j+t}-\sum_{u=1}^{j-i}b_{i+u}=x(\La',j'-i).\hspace{3.65cm}
$$
\end{proof}

\section{The case of at most two elements}\label{Kcase2}

In this section we will prove the Eisenbud-Goto conjecture in case that one of the equivalence classes which determine the regularity of $K[B]$ has at most two elements. As usual we define $\deg m=\sum_{j=1}^db_j$ for a monomial $m=y_1^{b_1}\cdot\ldots\cdot y_d^{b_d}$ in $T$.

\begin{defi}

We define the set $\Gamma(B)\subseteq\{\Gamma_1,\ldots,\Gamma_f\}$ by $\Gamma_t\in\Gamma(B)$ for $t\in\{1,\ldots,f\}$ if $\reg K[B] = \reg I_t+\deg h_t$.

\end{defi}

We have $\Gamma(B)\not=\emptyset$ by Equation~(\ref{Kregberechnung}). We note that the ideals and shifts which correspond to the elements of $\Gamma(B)$ are computed by the function \texttt{regularityMA} in \cite{KBEN}.

\begin{satz}\label{Keg2}

Let $\Gamma_t\in\Gamma(B)$ for some $t\in\{1,\ldots,f\}$. If $\#\Gamma_t\leq2$, then
$$
\reg K[B]\leq \deg K[B]-\codim K[B].
$$

\end{satz}
\begin{proof}

By construction we need to show that $\reg I_t + \deg h_t\leq\deg K[B]-c$. If $\#\Gamma_t=1$, then the assertion follows from Proposition~\ref{Kegred}. We therefore may assume that $\#\Gamma_t=2$, say, $\Gamma_t=\{x,x'\}$. We set $m:=y^{(x-h_t)/\alpha}$ and $n:=y^{(x'-h_t)/\alpha}$, that is, $I_t={(m,n)T}$. Since $I_t$ is minimally generated by $m$ and $n$ (see Remark \ref{KmingensK}) and by construction of $h_t$ we obtain that $m,n$ is a regular sequence on $T$. Using the Koszul Complex of $m,n$ (see, for example, \cite[Section~A2F]{KTGOS}) we get that the minimal graded free resolution of $I_t$ has the form
$$
0\longrightarrow T(-(\deg m+\deg n))\longrightarrow T(-\deg m) \oplus T(-\deg n)\longrightarrow I_t\longrightarrow 0.
$$
Then by \cite[Proposition]{KEG}
\begin{equation}\label{Keqa}
\reg K[B]= \reg I_t + \deg h_t=\deg m+\deg n-1+\deg h_t= \deg x+\deg x'-\deg h_t -1.
\end{equation}
Let $\lambda\in\Lambda_x$ and $\nu\in\LA_{x'}$. Note that $\La$ and $\n$ are crossless by Corollary~\ref{Kcross2el}. We set
$$
L:=\{x(\lambda,0),\ldots,x(\La,\deg x-2),x(\La,\deg x)\}\cup\{x'(\n,0),\ldots,x'(\n,\deg x'-2),x'(\n,\deg x')\}.
$$
We have $L\subseteq B_A$ by Lemma~\ref{K*prop}~$1)$. By construction we get that every element of $L$ is not equivalent to an element of $\{a_1,\ldots,a_c\}$, since for all $z\in L$ we have $\deg z\not=1$; see Remark~\ref{Ka1bisc}. Moreover, note that $x(\La,\deg x-1),x'(\n,\deg x'-1)\in\{a_1,\ldots,a_c\}$. Denote by $g$ the number of equivalence classes on $L$. We obtain
\begin{equation}\label{Keqaa}
g\stackrel{\mbox{\scriptsize{\ref{K*prop}}}}{=}\deg x+\deg x'-\#\left(\Delta(\La,\n)\setminus \{(\deg x-1,\deg x'-1)\}\right)\hspace{3.4cm}
\end{equation}
$$
\hspace{1.87cm}\geq \deg x+\deg x' - \#\Delta(\lambda,\n) =\deg x+\deg x'-\delta(\lambda,\n)-2\stackrel{\mbox{\scriptsize{\ref{Khungleichung}}}}{\geq} \deg x+\deg x' -\deg h_t -1,
$$
since $h(x,x')=h_t$. We have $\deg K[B]=f$ by Equation~(\ref{Kgleichung1}), and hence
$$
\deg K[B]\stackrel{\mbox{\scriptsize{\ref{Ka1bisc}}}}{\geq} g+c\stackrel{\mbox{\scriptsize{(\ref{Keqaa})}}}{\geq} \deg x+\deg x' -\deg h_t -1+c\stackrel{\mbox{\scriptsize{(\ref{Keqa})}}}{=}\reg K[B]+c.
$$
\end{proof}

We therefore obtain from Theorem~\ref{Keg2}:

\begin{fol}\label{Kcase22}

If $\#\Gamma_t\leq 2$ for all $t=1,\ldots,f$, then
$$
\reg K[B]\leq \deg K[B]-\codim K[B].
$$

\end{fol}

\begin{bsp}

Consider the following semigroup in $\mathbb N^4$ with $\alpha=6$:
$$
B=\erz{e_1,\ldots,e_4,(0,2,0,4), (3,0,2,1), (0,2,2,2)}.
$$
Using the function \texttt{regularityMA} in \cite{KBEN} we obtain $\reg K[B]=6$. Moreover, we have
$$
\Gamma_{t}=\{(3,6,4,11),(15,0,10,5)\}\in\Gamma(B),
$$
for some $t\in\{1,\ldots,f\}$, since $\reg I_t+\deg h_t=\reg {(y_2y_4,y_1^2y_3)T}+2=6$, hence the Eisenbud-Goto conjecture holds for $K[B]$ by Theorem~\ref{Keg2}. Note that there are also functions available testing the Buchsbaum, Cohen-Macaulay, Gorenstein, normal, and the seminormal property in the simplicial case. Using this we get that $K[B]$ is not Buchsbaum and not seminormal. 

\end{bsp}

\begin{bsp}

Let $\Gamma_t\in\Gamma(B)$ for some $t\in\{1,\ldots,f\}$ with $\#\Gamma_t>2$. Unfortunately this case is much more complicated. Consider the following possible equivalence class $\Gamma_t$ for $\alpha=20$:
$$
\Gamma_t=\{x=(44,104,12),y=(104,44,12),z=(24,24,72)\}.
$$
We get $h(x,y)=(44,44,12)$, $h(x,z)=(24,24,12)$, $h(y,z)=(24,24,12)$, and $h_t=(24,24,12)$. Assume that Conjecture~\ref{Kkomisch} holds, so $x$ and $y$ could have $4$ non-trivial pairwise equivalent elements, $x$ and $z$ could have $2$, as well as $y$ and $z$. Let us consider a worst case scenario:
$$
\begin{xy}
  \xymatrix{
x \ar@{~}[d]^{}\ar@{-}[r]^{} & x(\La,1) \ar@{-}[r]^{}  & x(\La,2) \ar@{-}[r]^{}\ar@{~}[d]^{}  & x(\La,3)  \ar@{-}[r]^{}  & x(\La,4)\ar@{~}[d]^{}   \ar@{-}[r]^{} & x(\La,5)\ar@{~}[d]^{}   \ar@{-}[r]^{} & x(\La,6)\ar@{~}[d]^{}  \ar@{-}[r]^{}  &  0 \ar@{~}[d]^{} \\
y \ar@{~}[d]^{}\ar@{-}[r]^{} & y(\n,1)\ar@{~}[d]^{}   \ar@{-}[r]^{} & y(\n,2) \ar@{-}[r]^{}  & y(\n,3)  \ar@{-}[r]^{} \ar@{~}[d]^{} & y(\n,4)  \ar@{-}[r]^{} & y(\n,5) \ar@{-}[r]^{}  & y(\n,6)  \ar@{-}[r]^{} & 0\ar@{~}[d]^{} \\
z \ar@{-}[r]^{} & z(\mu,1) \ar@{-}[r]^{}  & z(\mu,2) \ar@{-}[r]^{} \ar@/^ 0.0cm/@{~}[uul]^{}& z(\mu,3)  \ar@{-}[r]^{}  & z(\mu,4)\ar@/^ 0.0cm/@{~}[uul]^{}  \ar@{-}[rrr]^{}&&&  0
  }
\end{xy}
$$
for some $\La\in\LA_x, \n\in\LA_y$, and $\mu\in\LA_z$. Note that no element in the picture has degree $1$. If we follow the proof of Theorem~\ref{Keg2} we would get $g=10$. So $\reg I_t+\deg h_t$ should be smaller or equal to $10$. But this is not the case, since $\deg h_t=3$ and $\reg I_t=\reg{(y_1y_2^4,y_1^4y_2,y_3^3)T}=9$.

\end{bsp}

\section{Monomial curves}\label{Kmonomial}

In the following we will consider the case of monomial curves. So throughout this section we will assume that $d=2$, that is, $B$ is the submonoid of $(\mathbb N^2,+)$ which is generated by a set $\{e_1,e_2,a_1,\ldots,a_c\}\subset\mathbb N^2$, where $e_1=(\alpha,0)$, $e_2=(0,\alpha)$, ${a_i}_{[1]}+{a_i}_{[2]}=\alpha$ for all $i=1,\ldots,c$, and \mbox{$\gcd({a_1}_{[1]},\ldots,{a_c}_{[1]},\alpha)=1$}. It follows that $G(B)=\{(a,b)\in\mathbb Z^2\mid \alpha \mid (a+b)\}$ and hence
$$
f=\#(G(B)\cap D)=\#\{(a,b)\in\mathbb N^2 \mid a<\alpha, b<\alpha, \alpha\mid (a+b)\}=\alpha,
$$
that is, the number of equivalence classes on $B_A$ is equal to $\alpha$. We get $\deg K[B]=\alpha$ by Equation~(\ref{Kgleichung1}), moreover, $T=K[y_1,y_2]$. We can order the monomials of $T$ with respect to the lexicographic order, in particular, we can write every monomial ideal $I$ in $T$ in terms of its minimal monomial generators as follows:
$$
I = {(m_1,\ldots,m_r)T}, \mbox{with } m_i=y_1^{b_i} y_2^{c_i}, i=1,\ldots,r,
$$
where $b_1>\ldots >b_r\geq0$ and $0\leq c_1<\ldots< c_r$ (see \cite[Page~42]{KCCA}). We get $\reg I=\deg m_1$ if $r=1$. In case that $r\geq2$ it is well known that the regularity of $I$ can be computed by

\begin{prop}\label{Kregmon} 

Let $I$ be given as above and let $r\geq2$. We have
$$
\reg I =\max\limits_{i=1,\ldots,r-1}\maxk{b_i+c_{i+1}}-1.
$$

\end{prop}
\begin{proof}

Consider $g: T^r\twoheadrightarrow I$ where $\hat e_i\mapsto m_i$ for all $i$. By \cite[Proposition~3.1]{KCCA} the kernel of $g$ is minimally generated by $y_2^{c_{i+1}-c_i}\hat e_i-y_1^{b_i-b_{i+1}}\hat e_{i+1}$, $i=1,\ldots,r-1$, moreover, a minimal graded free resolution of $I$ has length $1$. Let $l\in\{1,\ldots,r-1\}$. Since $y_2^{c_{l+1}-c_l}\in T(-(b_l+c_l))_{b_l+c_{l+1}}$ and $y_1^{b_l-b_{l+1}} \in T(-(b_{l+1}+c_{l+1}))_{b_l+c_{l+1}}$ the minimal graded free resolution of $I$ has the form
$$
0\longrightarrow \bigoplus\limits_{j=1}^{r-1} T(-(b_j+c_{j+1})){\longrightarrow} \bigoplus\limits_{i=1}^{r} T(-(b_i+c_i)) {\longrightarrow} I\longrightarrow 0.
$$
We have $b_l+c_{l+1}-1 \geq \max\maxk{b_l+c_l, b_{l+1}+c_{l+1}}$ and therefore (see \cite[Proposition]{KEG})
$$
\reg I =\max\maxk{b_1+c_1,\ldots,b_r+c_r, b_1+c_2 - 1,\ldots, b_{r-1} + c_r -1}{=} \max\limits_{i=1,\ldots,r-1}\maxk{b_i+c_{i+1}}-1.
$$
\end{proof}

\begin{bem}\label{Klookideal}

Let $t\in\{1,\ldots,\alpha\}$ and $x,y\in\Gamma_t$. In case that $x\not=y$ there exist $i,j\in\{1,2\}$ with $x_{[i]}>y_{[i]}$ and $x_{[j]}<y_{[j]}$ by Remark~\ref{KmingensK}. Hence we can write the monomial ideal $I_t$ as
$$
I_t = {(m_1,\ldots,m_{\#\Gamma_t})T}, \mbox{with } m_i=y_1^{b_i} y_2^{c_i}\in\tilde\Gamma_t, i=1,\ldots,{\#\Gamma_t},
$$
where $b_1>\ldots >b_{\#\Gamma_t}\geq0$ and $0\leq c_1<\ldots< c_{\#\Gamma_t}$, moreover, by construction of $h_t$ we obtain $b_{\#\Gamma_t}=0$ and $c_1=0$. In case that $I_t$ is given as above we can easily construct $\Gamma_t$, more precisely, $\Gamma_t=\{(\alpha \cdot b_1,\alpha\cdot c_1)+h_t,\ldots,(\alpha \cdot b_{\#\Gamma_t},\alpha \cdot c_{\#\Gamma_t})+h_t\}$.

\end{bem}

In view of the staircase diagram of $I_t$ (see, for example, \cite[Page~42]{KCCA}) the monomials $m_i$ and $m_{i+1}$ are somehow adjacent. We translate this to the semigroup as follows:

\begin{defi}

Let $x,y\in\Gamma_t$ for some $t\in\{1,\ldots,\alpha\}$ with $x\not=y$, that is, $x_{[i]}>y_{[i]}$ and $x_{[j]}<y_{[j]}$ for some $i,j\in\{1,2\}$ with $i\not=j$. We say that $x$ and $y$ are \emph{adjacent} if there is no element $z\in\Gamma_t$ with $x_{[i]}>z_{[i]}>y_{[i]}$ and $x_{[j]}<z_{[j]}<y_{[j]}$.

\end{defi}

\begin{bsp}\label{Ksmooth5bsp}

Consider the following smooth monomial curve in $\mathbb P^5$ given by
$$
B=\erz{(12,0),(0,12),(11,1),(9,3),(4,8),(1,11)}.
$$
Then by \cite[Corollary~3.9]{KDAMNP} we obtain $\reg K[B]=4$. By using \cite{KBEN} we get
$$
K[B]\cong T\oplus T(-1)^4\oplus {(y_1,y_2)T}(-1)^2\oplus {(y_1,y_2^2)T}(-1)^2\oplus {(y_1^2,y_2)T}(-1)^2\oplus\underbrace{{(y_1^2,y_1y_2,y_2^3)T}}_{=I_{12}}(-1).
$$
Moreover, by Proposition~\ref{Kregmon} we have $\Gamma(B)=\{\Gamma_{12}\}$, where $\Gamma_{12}=\{(31,5),(19,17),(7,41)\}$. We note that $(31,5)$ and $(19,17)$ are adjacent, as well as $(19,17)$ and $(7,41)$.

\end{bsp}

\begin{bem}

Let us consider the case of smooth monomial curves, that is, we assume that $a_1=(\alpha-1,1)$ and $a_c=(1,\alpha-1)$. In this case we get that $\reg K[B]\leq \#L +1$ by \cite{KCMRMC}, where $\#L$ is the maximal number of consecutive integer points on the line $[(\alpha,0),(0,\alpha)]$ not belonging to $B$. This bound is clearly better than the L'vovsky bound \cite{KLV} (see introduction), anyway, even this bound is far from sharp, see \cite[Introduction]{KDAMNP}. Consider again Example~\ref{Ksmooth5bsp}, we get $\reg K[B]\leq 8$ by the Eisenbud-Goto bound, $\reg K[B]\leq 7$ by the L'vovsky bound, and $\reg K[B]\leq 5$ by the Hellus-Hoa-St\"uckrad bound. We will now give a short proof of the Eisenbud-Goto conjecture for smooth monomial curves:
 
Let $\Gamma_t\in\Gamma(B)$ for some $t\in\{1,\ldots,\alpha\}$. By Theorem~\ref{Keg2} (or Proposition~\ref{Kegred}) we may assume that $\#\Gamma_t\geq2$. Since $(\alpha-1,1),(1,\alpha-1)\in B$ we have $(k\alpha-l,l),(\alpha-l,k'\alpha+l)\in\Gamma_t$ for some $l,k,k'\in\mathbb N$ with $0<l<\alpha$. Set $x:=(k\alpha-l,l)$ and $x':=(\alpha-l,k'\alpha+l)$. Since $0<l<\alpha$ we obtain $h_t=(\alpha-l,l)$ and therefore $I_t={(y_1^{\deg x-1},\ldots,y_2^{\deg x'-1})T}$. We get
\begin{equation}\label{Ksmoothgleichung1}
\reg K[B]=\reg I_t+\deg h_t=\reg{(y_1^{\deg x-1},\ldots,y_2^{\deg x'-1})T}+1\stackrel{\mbox{\scriptsize{\ref{Kregmon}}}}{\leq} \deg x+\deg x'-2.
\end{equation}
Let $\Gamma_1=\{0\}$; note that $\deg h_{t'}=1$ for all $t'=2,\ldots,\alpha$ by the above argument. Fix $\La\in\LA_x$ and $\n\in\LA_{x'}$. Suppose that $x(\La,m)\sim x'(\n, n)$ for some $m\in\{1,\ldots,\deg x-1\}$ and some $n\in\{1,\ldots,\deg x'-1\}$. It follows that $\deg h(x,x')\geq 2$, by Lemma~\ref{K*prop}~1), Lemma~\ref{Klemeqcross}~$1)$, Equation~(\ref{K*prop2}), and Remark~\ref{Kungleichungh}, since $\deg h(z,z')\geq 1$ for all $z,z'\in B_A\setminus\{0\}$ with $z\sim z'$. This contradicts $h(x,x')=(\alpha-l,l)$ and hence $\#\Delta(\La,\n)=2$ (that is, Conjecture~\ref{Kkomisch} holds for $x$ and $x'$). By a similar argument as in Theorem~\ref{Keg2} we get
$$
\deg K[B]\geq \deg x + \deg x' - \#\Delta(\La,\n) + c \stackrel{\mbox{\scriptsize{(\ref{Ksmoothgleichung1})}}}{\geq} \reg K[B]+c.
$$

Let us consider the smooth rational monomial curves in $\mathbb{P}^{3}$, that is, we consider semigroups of the form $B=\erz{(\alpha,0),(0,\alpha),(\alpha-1,1),(1,\alpha-1)}$. We have $(\alpha-1,1)+(1,\alpha-1)\notin B_A$, hence $B_A=\{i(1,\alpha-1),j(\alpha-1,1)\}\mid 0\leq i,j\leq\alpha-2\}$. We obtain
$$
\Gamma_1=\{0\},\Gamma_2=\{(1,\alpha-1)\},\Gamma_{3}=\{(\alpha-1,1)\},\Gamma_{4}=\{(2,2\alpha-2),((\alpha-3)\alpha+2,\alpha-2)\},
$$
$$
\Gamma_5=\{(3,3\alpha-3),((\alpha-4)\alpha+3,\alpha-3)\},\ldots,\Gamma_{\alpha}=\{(\alpha-2,(\alpha-3)\alpha+2),(2\alpha-2,2)\}.\hspace{0.12cm}
$$
Hence
$$
K[B]\cong T\oplus T(-1)^2\oplus {(y_1^{\alpha-3},y_2)T}(-1)\oplus {(y_1^{\alpha-4},y_2^2)T}(-1)\oplus\ldots\oplus {(y_1,y_2^{\alpha-3})T}(-1),
$$
meaning each $T$-module of the form ${(y_1^{\beta},y_2^{\gamma})T}(-1)$, $1\leq \beta,\gamma\leq\alpha-3$ with $\beta+\gamma=\alpha-2$ appears exactly once in the decomposition. In case that $\alpha\geq3$ we have $\reg K[B]=\alpha-2=\deg K[B]-\codim K[B]$, that is, the Eisenbud-Goto conjecture is sharp in this case.

\end{bem}

\pagebreak

\begin{prop}\label{Kcrosslessbar}

If Conjecture~\ref{Kkomisch} holds for adjacent elements, then
$$
\reg K[B]\leq \deg K[B]-\codim K[B].
$$

\end{prop}
\begin{proof}

Let $\Gamma_t\in\Gamma(B)$ for some $t\in\{1,\ldots,\alpha\}$, moreover, let $I_t={(m_1,\ldots,m_{\#\Gamma_t})T}$, where $m_i=y_1^{b_i}y_2^{c_i}\in\tilde\Gamma_t$, $i=1,\ldots,\#\Gamma_t$. By Theorem~\ref{Keg2} (or Proposition~\ref{Kegred}) we may assume that $\#\Gamma_t\geq2$, moreover, we may assume that $b_1>\ldots >b_{\#\Gamma_t}=0$ and $0= c_1<\ldots< c_{\#\Gamma_t}$, see Remark~\ref{Klookideal}. By Proposition~\ref{Kregmon} we can fix a number $k\in\{1,\ldots,{\#\Gamma_t}-1\}$ such that $\reg I_t=b_k+c_{k+1}-1$. Fix $x,x'\in\Gamma_t$ such that $m_k=y^{(x-h_t)/\alpha}$ and $m_{k+1}=y^{(x'-h_t)/\alpha}$. Since the minimal monomial generators of $I_t$ are ordered with respect to the lexicographic order we get that $x$ and $x'$ are adjacent, see again Remark~\ref{Klookideal}. We have 
\begin{equation}\label{Kgleich11} 
\reg K[B]=\reg I_t +\deg h_t = b_k+c_{k+1}-1+\deg h_t\hspace{5.17cm}
\end{equation}
$$
\hspace{1.76cm}= ((x-h_t)/\alpha)_{[1]}+((x'-h_t)/\alpha)_{[2]}-1+\deg h_t=\deg {(x_{[1]}, x'_{[2]})}-1.
$$
Fix $\lambda\in\LA_x$ and $\n\in\LA_{x'}$ such that $\delta(x,x')=\delta(\lambda,\n)$, and consider the set
$$
L:=\{x(\lambda,0),\ldots,x(\La,\deg x-2),x(\La,\deg x)\}\cup\{x'(\n,0),\ldots,x'(\n,\deg x'-2),x'(\n,\deg x')\}.
$$
We have $L\subseteq B_A$ by Lemma~\ref{K*prop}~$1)$. By construction we get that every element of $L$ is not equivalent to an element of $\{a_1,\ldots,a_c\}$, since for all $z\in L$ we have $\deg z\not=1$; see Remark~\ref{Ka1bisc}. Moreover, note that $x(\La,\deg x-1),x'(\n,\deg x'-1)\in\{a_1,\ldots,a_c\}$. Denote by $g$ the number of equivalence classes on $L$. We obtain
\begin{equation}\label{Kgleich22}
g\stackrel{\mbox{\scriptsize{\ref{K*prop}}}}{=} \deg x+\deg x' - \#\left(\Delta(\La,\n)\setminus\{(\deg x-1,\deg x'-1)\}\right)\geq \deg x+\deg x'-\#\Delta(\lambda,\n)
\end{equation}
$$
\hspace{1.72cm}=\deg {(x_{[1]}, x'_{[2]})}+\deg {(x'_{[1]},x_{[2]})}-\delta(x,x')-2\stackrel{\mbox{\scriptsize{\ref{Kkomisch}}}}{\geq} \deg {(x_{[1]}, x'_{[2]})}-1,
$$
since $h(x,x')=(x'_{[1]},x_{[2]})$, and since $x$ and $x'$ are adjacent. We have $\deg K[B]=f=\alpha$, hence
$$
\deg K[B]\stackrel{\mbox{\scriptsize{\ref{Ka1bisc}}}}{\geq} g+c\stackrel{\mbox{\scriptsize{(\ref{Kgleich22})}}}{\geq} \deg {(x_{[1]}, x'_{[2]})}-1+c\stackrel{\mbox{\scriptsize{(\ref{Kgleich11})}}}{=} \reg K[B]+c.
$$
\end{proof}

To prove Conjecture~\ref{Kkomisch} for adjacent elements, we will now discuss the meaning of a cross.

\begin{bem}\label{Kcrosslook1}

Let $x,y\in B_A$ with $x\sim y$ and $x\not=y$, and assume that $x_{[1]}>y_{[1]}$ and $x_{[2]}<y_{[2]}$, see Remark~\ref{Klookideal}. Moreover, let $\La\in\LA_x$ and $\nu\in\LA_y$ be not crossless, that is,
$$
\begin{xy}
  \xymatrix{
x \ar@{-}[rr]^{} \ar@{~}[d] &  & x(\lambda,i)\ar@{-}[rr]^{}\ar@{~}[drr] & & x(\lambda,j)\ar@{-}[rr]^{} &  & 0 \ar@{~}[d]\\
y \ar@{-}[rr]^{} & & y(\nu,l)\ar@{-}[rr]^{}\ar@{~}[urr] & & y(\nu,k)  \ar@{-}[rr]^{}  & & \,0,
  }
\end{xy}
$$
for some $i,j,l,k\in\mathbb N$ with $0<i<j<\deg x$ and $0<l<k<\deg y$, see also Lemma~\ref{Kcrosslem}. 

Fix $i,k$ (we could also fix $l,j$), then we have one of the following cases:

\begin{enumerate}

\item[$1)$] $x(\La,i)_{[1]}>y(\n,k)_{[1]}$ and $x(\La,i)_{[2]}<y(\n,k)_{[2]}$,

\item[$2)$] $x(\La,i)_{[1]}<y(\n,k)_{[1]}$ and $x(\La,i)_{[2]}>y(\n,k)_{[2]}$,

\end{enumerate}
by Lemma~\ref{Kcrosslem}. The first case is what you normally would expect, since $x_{[1]}>y_{[1]}$. The second case looks a little strange, but still possible. Keep in mind that we have $x(\La^*,\deg x-i)\sim y(\n^*,\deg y-k)$ by Lemma~\ref{Klemeqcross}~$1)$, $x(\La^*,\deg x-i),y(\n^*,\deg y-k)\in B_A$ by Lemma~\ref{K*prop}~$1)$,  and $x(\La^*,\deg x-i)\not= y(\n^*,\deg y-k)$ by assertion $2)$ and $3)$ of Lemma~\ref{Kcrosslem}. Moreover, we have $x(\La^*,\deg x-i)+x(\La,i)=x$ and $y(\n^*,\deg y-k)+y(\n,k)=y$; see Equation~(\ref{K*prop2}).

\end{bem}

\begin{lem}\label{Kausklammern}

Consider the same situation as in Remark~\ref{Kcrosslook1}. Moreover, let $x$ and $y$ be adjacent. If $x(\La,i)_{[1]}>y(\n,k)_{[1]}$ and $x(\La,i)_{[2]}<y(\n,k)_{[2]}$, then
$$
x(\La^*,\deg x-i)_{[1]}<y(\n^*,\deg y-k)_{[1]} \,\mbox{ and }\, x(\La^*,\deg x-i)_{[2]}>y(\n^*,\deg y-k)_{[2]}.
$$

\end{lem}
\begin{proof}

Suppose that the assertion is not true. Applying Lemma~\ref{Kcrosslem}~$4)$ and $5)$ to the second assertion in the same lemma we get $x(\La^*,\deg x-i)_{[1]} > y(\n^*,\deg y-k)_{[1]}$ as well as $x(\La^*,\deg x-i)_{[2]}<y(\n^*,\deg y-k)_{[2]}$. We set $z:=y(\n,k)+x(\La^*,\deg x-i)$; by construction $z\in B$ and $z\sim x,y$, see also Remark~\ref{Kcrosslook1}. By Equation~(\ref{K*prop2}) we obtain
$$
x_{[1]}>z_{[1]}>y_{[1]} \,\mbox{ and }\, x_{[2]}<z_{[2]}<y_{[2]}.
$$
Consider an element $z':=z-n_1e_1-n_2e_2\in B$ for some $n_1,n_2\in\mathbb N$ such that $n_1+n_2$ is maximal. We have $z'\in B_A$, $z'\leq z$, $z'\not=x,y$, and $z'\sim z\sim x,y$. Suppose $z'_{[1]}\leq y_{[1]}$ or $z'_{[2]}\leq x_{[2]}$, then $z'\leq y$ or $z'\leq x$, which contradicts $y\in B_A$ or $x\in B_A$. Hence
$$
x_{[1]}>z'_{[1]}>y_{[1]} \,\mbox{ and }\, x_{[2]}<z'_{[2]}<y_{[2]},
$$
thus, $x$ and $y$ are not adjacent, which contradicts our assumption.
\end{proof}

We will now prove Conjecture~\ref{Kkomisch} for adjacent elements. This together with Proposition~\ref{Kcrosslessbar} confirms the Eisenbud-Goto conjecture for arbitrary monomial curves.

\begin{prop}\label{Kkomischclose}

Let $x,y\in B_A$ with $x\sim y$ and $x\not=y$. If $x$ and $y$ are adjacent, then
$$
\delta(x,y)\leq \deg h(x,y)-1,
$$
that is, Conjecture~\ref{Kkomisch} holds for adjacent elements.

\end{prop}
\begin{proof}

By Remark~\ref{Klookideal} we may assume that $x_{[1]}>y_{[1]}$ and $x_{[2]}<y_{[2]}$. Moreover, we may assume that $x$ and $y$ are not crossless by Corollary~\ref{Kcrosslesskor}. Let us fix a maximal cross in the following sense: let $(\La,\n,i,j,l,k)\in\LA_x\times\LA_y\times \mathbb N^4$ be a cross such that $j-i$ is maximal among all crosses; say, $\La=(b_1,\ldots,b_{\deg x})$ and $\n=(g_1,\ldots,g_{\deg y})$. This can be illustrated by the following picture:
$$
\begin{xy}
  \xymatrix{
x \ar@{-}[rr]^{} \ar@{~}[d] &  & x(\lambda,i)\ar@{-}[rr]^{}\ar@{~}[drr] & & x(\lambda,j)\ar@{-}[rr]^{} &  & 0 \ar@{~}[d]\\
y \ar@{-}[rr]^{} & & y(\nu,l)\ar@{-}[rr]^{}\ar@{~}[urr] & & y(\nu,k)  \ar@{-}[rr]^{}  & & \,0.
  }
\end{xy}
$$
Keep in mind that if $(m,n)\in\Delta(\La,\nu)$, then $(m,n')\notin\Delta(\La,\nu)$ for all $n'\in\{0,\ldots,\deg y\}\setminus\{n\}$, see Remark~\ref{Kprojektion}. Moreover, recall that every permutation of a sequence with $*$-property is again a sequence with $*$-property by Remark~\ref{Kpermutation}. Without loss of generality we therefore may assume that for all $j',k'\in\mathbb N$ with $j<j'<\deg x$ and $k<k'<\deg y$ we have $x(\La, j')\not\sim y(\n,k')$; otherwise we can replace $\La$ by $\La'$ and $\nu$ by $\nu'$, where
$$
\La'=(b_{j'+1},\ldots,b_{\deg x},b_1,\ldots,b_{j'})\in\LA_x\quad \mbox{and}\quad\n'=(g_{k'+1},\ldots,g_{\deg y},g_1,\ldots,g_{k'})\in\LA_y,
$$
since $(\La',\n',\deg x-j'+i,\deg x-j'+j,\deg y-k'+l,\deg y-k'+k)$ is again a cross of height $(j-i,k-l)$, and since such a replacement can only be done finitely many times. Let $i'\in\mathbb N$ with $i'<i$ be maximal such that $x(\La,i')\sim y(\n,l')$ for some $l'\in\{0,\ldots,\deg y\}$. Set $x':=x(\La^*,\deg x-i')$, $y':=y(\n^*,\deg y-l')$, $x'':=x(\La,i')$, $y'':=y(\n,l')$, $\La':=(b_1,\ldots,b_{i'})\in\LA_{x'}$, and $\n':=(g_1,\ldots,g_{l'})\in\LA_{y'}$ (see Remark~\ref{Kaussehen1}). In case that $i'=0$ we set $\#\Delta(\La',\n'):=1$, note that $i'=0$ if and only if $l'=0$ by Lemma~\ref{K*prop}~$2)$. We have $x'\sim y'$ by Lemma~\ref{Klemeqcross}~$1)$, and $x=x'+x''$ as well as $y=y'+y''$ by Equation~(\ref{K*prop2}).

We claim that:
\begin{equation}\label{K1224}
\#\Delta(\La,\n)\leq \#\Delta(\La',\n') + \deg(x(\La,i)-x(\La,j))+2.
\end{equation}
Consider an element $j'\in\mathbb N$ with $j<j'<\deg x$ and suppose that $x(\La,j')\sim y(\n,k')$ for some $k'\in\{0,\ldots,\deg y\}$. By construction we have $k'<k$. Hence we get a cross $(\La,\n,i,j',k',k)$ of height $(j'-i,k-k')$, which contradicts the maximality of $j-i$, and therefore $x(\La,j')\not\sim y(\n,k')$ for all $k'\in\{0,\ldots,\deg y\}$. We have (see Remark~\ref{Kprojektion})
$$
\#\Delta(\La,\n)\leq \#\left(\Delta(\La,\n)\cap(\{0,\ldots,i'\}\times\mathbb N)\right) + \deg(x(\La,i)-x(\La,j))+2,
$$
that is, we need to show that $\left(\Delta(\La,\n)\cap(\{0,\ldots,i'\}\times\mathbb N)\right)\subseteq \Delta(\La',\n')$. In case that $i'=0$ we have $\#\left(\Delta(\La,\n)\cap(\{0,\ldots,i'\}\times\mathbb N)\right)=1$, that is, we may assume that $i'>0$. Suppose that $l'>l$, hence $(\La,\n,i',j,l,l')$ is a cross of height $(j-i',l'-l)$, which contradicts the maximality of $j-i$. That means $l'<l$, since $l'\not=l$, that is,
$$
\begin{xy}
  \xymatrix{
x \ar@{~}[d]\ar@{-}[rrr]^{}&&& x(\La,i')\ar@{~}[d]\ar@{-}[r]^{} & x(\lambda,i)\ar@{-}[rr]^{}\ar@{~}[drr] & & x(\lambda,j)\ar@{-}[rr]^{} &  & 0 \ar@{~}[d]\\
y\ar@{-}[rrr]^{} & & & y(\n,l')\ar@{-}[r]^{} & y(\nu,l)\ar@{-}[rr]^{}\ar@{~}[urr] & & y(\nu,k)  \ar@{-}[rr]^{}  & & \,0.
  }
\end{xy}
$$
Let $(m,n)\in\left(\Delta(\La,\n)\cap(\{0,\ldots,i'\}\times\mathbb N)\right)$. Since $(0,0),(i',l')\in\left(\Delta(\La,\n)\cap(\{0,\ldots,i'\}\times\mathbb N)\right)$ and $(0,0),(i',l')\in\Delta(\La',\nu')$ we may assume that $m\notin\{0,i'\}$. Suppose that $x(\La,m)\sim y(\n,n)$ for some $n\geq l'$. By a similar argument as above we obtain $n<l$ and clearly $n\not=l'$, that is, we suppose that $l'<n<l$. Hence $(\La,\n,m,i',l',n)$ and $(\La,\n,i,j,l,k)$ are two crosses with $i'<i$ and $n<l$ which contradicts the maximality of $j-i$, see Proposition~\ref{Kgluingcrosses}. That means $n<l'$ and therefore $(m,n)\in\Delta(\La',\n')$ by Lemma~\ref{Klemeqcross}~$2)$, which proves Equation~(\ref{K1224}).

Since $x$ and $y$ are adjacent we get by Remark~\ref{Kcrosslook1} and Lemma~\ref{Kausklammern} one of the following cases

\begin{enumerate}

\item[$1)$] $x(\La^*,\deg x-i)_{[1]}<y(\n^*,\deg y-k)_{[1]}$ and $x(\La^*,\deg x-i)_{[2]}>y(\n^*,\deg y-k)_{[2]}$,

\item[$2)$] $x(\La,i)_{[1]}<y(\n,k)_{[1]}$ and $x(\La,i)_{[2]}>y(\n,k)_{[2]}$.

\end{enumerate}
\textbf{Case 1:}\\
Applying Lemma~\ref{Kcrossgeq} to Lemma~\ref{Kcrosslem}~$2)$ we get $x(\La^*,\deg x-j)_{[2]}>y(\n^*,\deg y-l)_{[2]}$. Moreover, by applying Lemma~\ref{Kcrosslem}~$4)$ to the second assertion in the same lemma it follows that $x(\La^*,\deg x-j)_{[1]}<y(\n^*,\deg y-l)_{[1]}$. By construction $h(x,y)=(y_{[1]},x_{[2]})$. Hence
$$
h(x,y)_{[1]}=y_{[1]}\geq y(\n^*,\deg y-l)_{[1]}>x(\La^*,\deg x-j)_{[1]}
$$
and 
$$
h(x,y)_{[2]}=x_{[2]}\geq x(\La^*,\deg x-j)_{[2]}.
$$
Thus,
\begin{equation}\label{K9877}
\deg h(x,y)-1 \geq \deg x(\La^*,\deg x-j).
\end{equation}
Moreover, we have $\Delta(\La',\n')\subseteq(\{0,\ldots,i'\}\times\{0,\ldots,l'\})$, that is, $\#\Delta(\La',\n')\leq i'+1$ (see Remark~\ref{Kprojektion}) and $i'+1\leq i$. We get
\begin{equation}\label{K999}
\#\Delta(\La',\n') + \deg(x(\La,i)-x(\La,j))\leq i + \deg x-i-(\deg x-j) = \deg x(\La^*,\deg x-j),
\end{equation}
and therefore
$$
\delta(x,y)\leq\delta(\La,\n)=\#\Delta(\La,\n)-2\stackrel{\mbox{\scriptsize{(\ref{K1224})}}}\leq\#\Delta(\La',\n') + \deg(x(\La,i)-x(\La,j))\stackrel{\mbox{\scriptsize{(\ref{K9877}),(\ref{K999})}}}{\leq} \deg h(x,y)-1.
$$
\textbf{Case 2:}\\
By Lemma~\ref{Klemeqcross}~$2)$ and Proposition~\ref{Kgluingcrosses} $\La'$ and $\n'$ are crossless, since $(j-i)$ is assumed to be maximal. Hence by Proposition~\ref{Khungleichung} we get:
\begin{equation}\label{K1223}
\#\Delta(\La',\n')-2\leq\deg h(x',y')-1.
\end{equation}
In case that $i'=0$ we have $\#\Delta(\La',\n')=1$ and $\deg h(x',y')=0$, that is, Equation~(\ref{K1223}) holds. We get $x''_{[2]}\geq x(\La,i)_{[2]}$, and $y''_{[1]}\geq y(\n,k)_{[1]}>x(\La,i)_{[1]}$ and therefore we obtain 
\begin{equation}\label{Kgleichungeing}
\deg{(y''_{[1]},x''_{[2]})}\geq\deg x(\La,i)+1.
\end{equation}
Hence
$$
\deg h(x,y)-1=\deg{(y_{[1]},x_{[2]})}-1=\deg{(y'_{[1]},x'_{[2]})}+\deg{(y''_{[1]},x''_{[2]})}-1\hspace{2.43cm}
$$
$$
\geq \deg h(x',y')-1+\deg{(y''_{[1]},x''_{[2]})}\stackrel{\mbox{\scriptsize{(\ref{K1223}),(\ref{Kgleichungeing})}}}{\geq}\#\Delta(\La',\n')-2+\deg x(\La,i)+1\hspace{1.83cm}
$$
$$
\hspace{0.15cm}\geq \#\Delta(\La',\n')-2+\deg(x(\La,i)-x(\La,j))+1+1\stackrel{\mbox{\scriptsize{(\ref{K1224})}}}{\geq} \#\Delta(\La,\n)-2=\delta(\La,\n)\geq \delta(x,y).
$$
\end{proof}

We therefore obtain from Proposition~\ref{Kcrosslessbar} and Proposition~\ref{Kkomischclose}:

\begin{satz}\label{Kegmonomial}

Let $\dim K[B]=2$. We have
$$
\reg K[B]\leq \deg K[B]-\codim K[B].
$$

\end{satz}

\begin{bem}

In view of our approach it would be nice to confirm Conjecture~\ref{Kkomisch}, or even better, the stronger version which is mentioned after the conjecture. Moreover, it would be very interesting to find a L'vovsky version of Proposition~\ref{Kegred}, meaning, we should find a good bound for the degree of an element of $B_A$ in terms of the two biggest gaps of $B$. This could help to find a combinatorial proof of the L'vovsky bound.

\end{bem}

\section*{Acknowledgement}

The basic ideas have been developed, while the author was visiting the University of Barcelona in December 2010; the first version of the paper was completed during a visit at the University of California and at the University of Utah in spring 2011. The author would like to thank Santiago Zarzuela, David Eisenbud, and Anurag Singh for their hospitality and for useful conversations. Moreover, he would like to thank his PhD advisor J\"urgen St\"uckrad for many helpful discussions as well as for suggesting working on simplicial affine semigroup rings and for supporting his travel activities. Finally, the author would like to thank the Max Planck Institute for Mathematics in the Sciences in Leipzig and the International Max Planck Research School Mathematics in the Sciences for their financial support.

\end{document}